\documentclass[11pt]{amsart}

\usepackage{amscd,amsmath,amssymb}
\usepackage[mathscr]{eucal}
\usepackage[all]{xy}
\usepackage[charter]{mathdesign}

\makeatletter \@addtoreset{equation}{section}\makeatother

\DeclareMathAlphabet{\mathpzc}{OT1}{pzc}{m}{it}

\newtheorem{theorem}{Theorem}[section]
\newtheorem{lemma}[theorem]{Lemma}
\newtheorem{proposition}[theorem]{Proposition}
\newtheorem{corollary}[theorem]{Corollary}

\newtheorem{remark}[theorem]{Remark}
\newtheorem{conjecture}[theorem]{Conjecture}
\newtheorem{example}[theorem]{Example}

\newcommand{\Z}{{\mathbb{Z}}}
\newcommand{\ZZ}{{\mathbb{Z}/2}}
\newcommand{\C}{{\mathbb{C}}}

\newcommand{\K}{{\C}}
\newcommand{\lm}{{\mu}}
\newcommand{\ld}{{\delta}}
\newcommand{\bd}{{b(\delta)}}
\newcommand{\bm}{{b(\mu)}}
\newcommand{\rC}{{\mathsf{C}}}
\newcommand{\rc}{{v}}

\newcommand{\bft}{{\underline{c}}}
\newcommand{\bfg}{{\underline{g}}}
\newcommand{\rCn}{{{\mathsf{C}}}}
\newcommand{\rCp}{{\mathsf{C}^\Pi_*}}
\newcommand{\rHH}{{\mathsf{HH}}}
\newcommand{\rHP}{{\mathsf{HP}}}
\newcommand{\rHN}{{\mathsf{HN}}}

\newcommand{\U}{{\mathnormal{e}(\delta)}}
\newcommand{\V}{{\mathnormal{E}(\delta)}}
\newcommand{\ck}{{\mathrm{sh}}}
\newcommand{\cK}{{\mathrm{Sh}}}

\newcommand{\cD}{{\mathscr{D}}}
\newcommand{\MF}{{\mathscr{M\!F}(f)}}

\newcommand{\sh}{{\mathtt{sh}}}
\newcommand{\Sh}{{\mathtt{Sh}}}
\newcommand{\End}{{\mathrm{End}}}

\newcommand{\STR}{{{\Upsilon}}}

\newcommand{\rH}{{\mathrm{H}}}
\newcommand{\cH}{{\mathsf{H}}}
\newcommand{\cA}{{\mathcal{A}}}

\newcommand{\uu}{{(\!(u)\!)}}
\newcommand{\uuu}{{[\![u]\!]}}

\newcommand{\str}{{\mathrm{str}}}

\newcommand{\met}{{{K}}}

\newcommand{\cI}{{\mathcal{I}_f}}
\newcommand{\bI}{{I_f}}

\newcommand{\Aw}{{A_f}}
\newcommand{\Dw}{{D_f}}
\newcommand{\Kw}{{K_f}}

\newcommand{\zz}{{[\![{z}_1,\ldots,{z}_k]\!]}}
\newcommand{\xx}{{[\![x_1,\ldots,x_n]\!]}}
\newcommand{\zzu}{{[\![{z}_1,\ldots,{z}_k,u]\!]}}
\newcommand{\conn}{{\boldsymbol{\nabla}}}
\newcommand{\bconn}{{\bigtriangledown}}

\allowdisplaybreaks

\title{Matrix factorizations and higher residue pairings}
\author{D. Shklyarov}
\address{Freiburg Institute for Advanced Studies (FRIAS) and Mathematics Institute, University of Freiburg,  Germany}
\email{dmytro.shklyarov@math.uni-freiburg.de}

\begin{document}
\begin{abstract}
The periodic cyclic homology of any proper dg category comes equipped with a canonical pairing. We show that in the case of the dg category of matrix factorizations of an isolated singularity the canonical pairing can be identified with Kyoji Saito's higher residue pairing on the twisted de Rham cohomology of the singularity. 
\end{abstract}

\maketitle

\baselineskip 1.6pc

\section{Introduction} 
\subsection{Higher residue pairings}
About thirty years ago  Kyoji Saito constructed first examples of what was later titled Frobenius manifolds.  An abstract Frobenius manifold is a manifold whose holomorphic tangent bundle is equipped with a fiberwise Frobenius algebra structure satisfying some flatness conditions. The underlying manifolds in K. Saito's examples are miniversal unfolding spaces of function germs. If $f$ is a germ of an analytic function at the origin ${\sf 0}\in\C^n$ such that ${\sf 0}$ is its isolated critical point and $f({\sf 0})=0$ then the tangent space (at the base-point) to the unfolding space of $f$ can be identified with $\C\{x_1,\ldots, x_n\}/(\partial_1f,\ldots, \partial_nf)$, the Milnor algebra of $f$. The latter is, in addition, a Frobenius algebra with respect to the residue pairing
\begin{equation}\label{res}
\eta_f(\xi_1, \xi_2)=\frac{1}{(2\pi i)^n}\oint\frac{\,\xi_1\,\xi_2\,d\underline{x}}{\partial_1f\ldots\partial_nf}, \quad d\underline{x}:=dx_1\wedge\ldots\wedge dx_n
\end{equation}
Roughly, the Frobenius manifold is constructed by ``spreading''  this Frobenius algebra structure over the space of deformations of $f$. The actual procedure involves a generalization of variations of Hodge structures and period maps. The reader is referred to \cite{ST} and \cite{Her} for a brief and exhaustive account, respectively. 

Let us recall some of the Hodge theoretic definitions relevant for the above theory.

The role of the de Rham cohomology is played in this setting by the twisted de Rham cohomology\footnote{Since $f$ has an isolated singularity, only the $n$-th twisted de Rham cohomology is non-trivial.}
\begin{equation}\label{tw}
\cH_f:=\rH^n(\Omega^{\ast}_{\C^n,{\sf 0}}\uu, -d{f} +ud)\,
\end{equation}
where $u$ is a variable and $\Omega^{\ast}_{\C^n,{\sf 0}}$ stands for the germs of analytic forms at ${\sf 0}$. 
$\cH_f$ carries a Hodge-like structure which comprises the following data:  
the connection ${\bconn}^{f}$ induced by the operator
\begin{equation}\label{clascon}
\partial_u+\frac{f}{u^2}:\Omega^{\ast}_{\C^n,{\sf 0}}\uu\to\Omega^{\ast}_{\C^n,{\sf 0}}\uu
\end{equation}
and the free full-rank $\C\uuu$-submodule $\cH_f^{(0)}:=\rH^n(\Omega^{\ast}_{\C^n,{\sf 0}}\uuu, -d{f}+ud)$ (more precisely, the image in $\cH_f$ of the latter module under the natural embedding of the complexes).

The previous two pieces of data should be viewed as a generalization of the classical Hodge filtration. In addition, one has a generalization of polarizations. Their role is played by the so-called {\it higher residue pairings} \cite{Sa} or, more precisely, by their ``generating function''  \footnote{It is this generating function that we will call the higher residue pairing in the present paper.}
\[
\Kw: \cH_f\times\cH_f\to \C\uu.
\] 
The latter is a $\C\uu$-sesquilinear ${\bconn}^{f}$-flat pairing 
\begin{eqnarray*}
&&\Kw(\lambda(u)\,\cdot\,, \,\cdot\,)=\Kw(\,\cdot\,, \lambda(u)^\star\,\cdot\,)=\lambda(u)\Kw(\,\cdot\,, \,\cdot\,)\quad (\lambda(u)\in\K\uu, \,\,\lambda(u)^\star:=\lambda(-u)), \\
&&\Kw({\bconn}_{{\!\!{\partial_u}}}^{f}\,\cdot\,, \,\cdot\,)-\Kw(\,\cdot\,, {\bconn}_{{\!\!{\partial_u}}}^{f}\,\cdot\,)=\partial_u\Kw(\,\cdot\,, \,\cdot\,)
\end{eqnarray*}
 which extends the residue pairing (\ref{res}) in the following sense: the restriction of $\Kw$ to $\cH_f^{(0)}$ takes values in $u^n\K\uuu$ and satisfies
\begin{eqnarray}\label{p03}
\Kw(\,\cdot\,,\, \cdot\,)=u^n\eta_f(\,\,\cdot\,|_{u=0}, \,\,\cdot\,|_{u=0}) +O(u^{n+1}).
\end{eqnarray}
Here ``$|_{u=0}$'' stands for the projection onto $\cH_f^{(0)}/u\cH_f^{(0)}\simeq\Omega^n_{\C^n,{\sf 0}}/d{f}\wedge\Omega^{n-1}_{\C^n,{\sf 0}}$; the pairing (\ref{res}) is transferred to the latter space using the isomorphism  $\xi\mapsto\xi\,d\underline{x}$.

As shown in \cite{LLS}, $\Kw$ can be "lifted" to the complex $(\Omega^{\ast}_{\C^n,{\sf 0}}\uu, -d{f} +ud)$. Namely, there is an explicit (albeit quite complicated) quasi-isomorphism from this complex to a similar complex where the space of holomorphic forms is replaced by that of $C^\infty$-forms with compact supports. The cochain-level pairing is obtained by combining the quasi-isomorphism with the ordinary integral on the space of compactly supported forms.

\subsection{A categorical counterpart of the higher residue pairing and the main theorem} 

Let us explain now what the present work is about. In a sense, it is sequel to \cite{Shk2} and we will start by recalling the contents of that paper. 

Henceforth, $f$ will be viewed as an element of $R:=\C\xx$ rather than an analytic function germ. Accordingly, we will replace $\Omega^{\ast}_{\C^n,{\sf 0}}$ with the space $\widehat{\Omega}^\ast=\widehat{\Omega}^\ast_{\C^n,{\sf 0}}$ of forms with formal coefficients (as shown in \cite{Sch}, this does not affect the twisted de Rham cohomology).

Following ideas outlined in \cite{KKP,Kon,KS}, we showed in \cite{Shk2} that the triple $(\cH_f, {\bconn}^{f}, \cH_f^{(0)})$  can be reconstructed, up to a ``Tate twist'', from a differential $\ZZ$-graded (dg) category ${\MF}$ of {\it matrix factorizations} of $f$. The objects of this category, we recall, are super-vector bundles $E=E^\pm$ on ${\rm Spec}\,R$ endowed with an odd endomorphism $D$ satisfying $D^2=f\cdot{\rm id}_E$. (A detailed description of ${\MF}$ can be found in \cite{Dyck}.) The starting points in \cite{Shk2} were 
\begin{itemize}
\item An explicit isomorphism $\bI$, found in \cite{Seg}, between the Hochschild homology $\rHH_n({\MF})$ and the space $\widehat{\Omega}^n/d{f}\wedge\widehat{\Omega}^{n-1}$. The very same formula turns out to produce isomorphisms
\begin{equation}\label{qis}
\bI: \rHP_n({\MF})\simeq \cH_f, \quad  \rHN_n({\MF})\simeq \cH_f^{(0)}
\end{equation}
where $\rHP_*$ and $\rHN_*$ stand for the periodic and negative cyclic homology, respectively\footnote{Since ${\MF}$ is $\ZZ$-graded, all its homology groups are also only $\ZZ$-graded; that is, the subscript $n$ should be understood as ``$n\,\text{mod}\, 2$''. Note also that $\rHH_{n+1}({\MF})$, $\rHN_{n+1}({\MF})$, $\rHP_{n+1}({\MF})$ are all trivial.}.
\item Existence of a canonical connection ${\bconn}^{\cD}$ on the periodic cyclic homology $\rHP_*({\cD})$ of any dg category ${\cD}$. The definition of the connection was previously sketched in \cite{KKP,Kon,KS}.
\end{itemize}
The main result of \cite{Shk2} is the formula
\begin{equation*}
({\bconn}_{{\!\!{\partial_u}}}^{f}-\frac{n}{2u})\cdot\bI=\bI\cdot{\bconn}_{{\!\!{\partial_u}}}^{{\MF}}
\end{equation*}
(The extra term in the left-hand side is the aforementioned Tate twist.) 

Our goal in the present work is to give a categorical interpretation of the higher residue pairing (and thereby fill in the question mark in the table at the end of \cite{ST}). Namely, the periodic cyclic homology $\rHP_*(\cD)$ of any {\it proper} dg category $\cD$  comes equipped with a canonical $\C\uu$-sesquilinear pairing $\met_\cD$. Here ``proper'' means ${\rm Hom}_\cD(X,Y)\in \mathcal{C}_{\rm fd}$ for all $X,Y\in\cD$ where $\mathcal{C}_{\rm fd}$ stands for the dg category of $\ZZ$-graded complexes of vector spaces with finite-dimensional cohomology \cite{KS}. The pairing is defined as the composition (cf. \cite{Shk1})
\begin{equation}\label{canpr}
\rHP_\ast(\cD)\otimes\rHP_\ast(\cD){\to}\rHP_\ast(\cD)\otimes\rHP_\ast(\cD^{\rm op}){\to}\rHP_\ast(\cD\otimes \cD^{\rm op}){\to}\rHP_\ast(\mathcal{C}_{\rm fd})\simeq\C\uu
\end{equation}
where the first map comes from a canonical $\C\uu$-antilinear isomorphism $\rHP_\ast(\cD){\to}\rHP_\ast(\cD^{\rm op})$, the second map is a well-known K\"unneth type isomorphism, and the third map is induced by the dg functor 
\[
\cD\otimes \cD^{\rm op}{\to}\mathcal{C}_{\rm fd},\quad (X,Y)\mapsto {\rm Hom}_\cD(Y,X).
\]
The isomorphism $\rHP_\ast(\mathcal{C}_{\rm fd})\simeq\C\uu$ is a consequence of the dg Morita invariance of the cyclic homology \cite{K2,K3}. (To be more precise, we need to fix a specific isomorphism: let us fix the one induced by the canonical map $\rHP_\ast(\C)\to \rHP_\ast(\mathcal{C}_{\rm fd})$.)

The category ${\MF}$ is known to be proper \cite{Dyck}. Our main result is  
\begin{theorem}\label{MC} There exists $\mathsf{const}\in\C$ such that 
\begin{equation*}
\Kw(\,\bI(\cdot)\,,\,\bI(\cdot)\,)=\mathsf{const}\cdot u^n\met_{\MF}(\,\cdot\,,\,\cdot\,).
\end{equation*}
\end{theorem}

Let us point out one corollary. Observe that the initial term  ``$\eta_\cD:=\met_\cD$ mod $u$'' of the canonical pairing is a well-defined pairing on the Hochschild homology $\rHH_\ast(\cD)$. (This pairing was studied in detail in \cite{Shk1}; to define it formally, replace $\rHP_\ast$ with $\rHH_\ast$ and $\C\uu$ with $\C$ in (\ref{canpr}).)
Then Theorem \ref{MC} and (\ref{p03}) imply ($\mathsf{const}$ below is the same as in the theorem)
\begin{corollary}\label{corol}
$\eta_f(\,\bI(\cdot)\,, \,\bI(\cdot)\,)=\mathsf{const}\cdot \eta_{\MF}(\,\cdot\,,\,\cdot\,)$.
\end{corollary}
Another (quite different) way to identify the canonical pairing on the Hochschild homology of the category $\MF$ with the pairing $\eta_f$ was found earlier by A.~Polishchuk and A.~Vaintrob; cf. \cite[Corollary 4.1.3]{PV}. 

We believe that the constant $\mathsf{const}$ is essentially independent of $f$ and can be computed by studying the case $f=x_1^2+x_2^2+\ldots+x_n^2$ which leads to   

\begin{conjecture}\label{cj} $\mathsf{const}=(-1)^{\frac{n(n+1)}{2}}$.
\end{conjecture}

\begin{remark}{\rm It would be interesting to find an explicit relationship between Corollary \ref{corol} and the aforementioned result of Polishchuk and Vaintrob. This should lead to the proof of the conjecture. Note that the analogous constant in \cite{PV} is slightly different, namely, it equals $(-1)^{\frac{n(n-1)}{2}}$. Most probably, this discrepancy is due to the fact that the canonical pairing, considered in {\it loc. cit.}, lives on $\rHH_\ast(\MF)\otimes\rHH_\ast(\MF^{\rm op})$ whereas our pairing lives on $\rHH_\ast(\MF)\otimes\rHH_\ast(\MF)$. There is an explicit canonical isomorphism $\rHH_\ast(\MF)\to\rHH_\ast(\MF^{\rm op})$ (see (\ref{Phi})), and we believe it is this isomorphism that relates the two constants.
}
\end{remark}

\subsection{On the proof of the main theorem} 
To begin with, an important technical aspect of the proof is the possibility to replace the category $\MF$ with its more economical model. Namely, as shown in \cite{Dyck}, the dg category ${\MF}$ is dg Morita equivalent to an explicit dg algebra $\Aw$ which, when viewed as a dg category with a single object, embeds into ${\MF}$. This embedding induces isomorphisms
\begin{equation}\label{ident}
\rHP_*(\Aw)\simeq \rHP_*(\MF), \quad\rHN_*(\Aw)\simeq \rHN_*(\MF),\quad \text{etc.}
\end{equation}
and these isomorphisms allow one to identify ${\bconn}^{{\MF}}$ with ${\bconn}^{{\Aw}}$ and $\met_{\MF}$ with $\met_{\Aw}$. In other words, ${\MF}$ and $\Aw$ share all the homological invariants we are interested in. Therefore, in the main body of the present paper we will be working in the setting of dg algebras rather than dg categories.
 
The proof of Theorem \ref{MC} is based on Morihiko Saito's uniqueness theorem for the higher residue pairings \cite[Appendix]{Sa1}. Let us recall what the theorem claims.

Fix a collection of elements $\bfg=\{g_i\}_{i=1,\ldots,k}$ in $\C\xx$ whose representatives form a basis in the Milnor algebra of $f$. Fix also some variables $z_1,\ldots, z_k$ and consider the free rank $k$ $\C\zzu$-module\footnote{It is easy to show that all the other cohomology groups of this complex are trivial, just as in the ``undeformed'' case.}
\begin{equation*}
\cH^{(0)}_{f,\bfg}:=\rH^n(\widehat{\Omega}^\ast\zzu, -df-z_1dg_1-\ldots-z_k dg_k+ud)
\end{equation*}
and its localization at $u$
\begin{equation*}
\cH_{f,\bfg}:=\rH^n(\widehat{\Omega}^\ast\zz\uu, -df-z_1dg_1-\ldots-z_k dg_k+ud).
\end{equation*}
The latter module carries the so-called (extended) Gauss-Manin connection ${\bconn}^{f,\bfg}$  induced by the operators 
\begin{eqnarray}\label{exGM}
{\partial_u}+\frac{f+z_1g_1+\ldots+z_k g_k}{u^2}, \,\,{\partial_{{z}_i}}-\frac{g_i}{u}: \,\,\widehat{\Omega}^\ast\zz\uu\to \widehat{\Omega}^\ast\zz\uu.
\end{eqnarray}
The higher residue pairing extends to a $\C\zz$-linear $\C\uuu$-sesquilinear $u^n\C\zzu$-valued pairing $K_{f,\bfg}$ on $\cH^{(0)}_{f,\bfg}$ in such a way that the induced pairing on $\cH_{f,\bfg}$ is ${\bconn}^{f,\bfg}\,\,$--flat:
\begin{eqnarray}\label{flat1}
\,  \\
K_{f,\bfg}({\bconn}_{{\!\!{\partial_{u}}}}^{f,\bfg}\,\,\cdot\,,\,\cdot\,)-K_{f,\bfg}(\,\cdot\,,{\bconn}_{{\!\!{\partial_{u}}}}^{f,\bfg}\,\,\cdot\,)={\partial_u}K_{f,\bfg}(\,\cdot\,,\,\cdot\,),\quad
K_{f,\bfg}({\bconn}_{{\!\!{\partial_{z_i}}}}^{f,\bfg}\,\,\cdot\,,\,\cdot\,)+K_{f,\bfg}(\,\cdot\,,{\bconn}_{{\!\!{\partial_{z_i}}}}^{f,\bfg}\,\,\cdot\,)={\partial_{z_i}}K_{f,\bfg}(\,\cdot\,,\,\cdot\,).\nonumber
\end{eqnarray}
The uniqueness theorem says that, up to a constant, $K_{f,\bfg}$ is the only pairing satisfying this flatness condition.

Let us briefly outline the proof of Theorem \ref{MC}.

The above elements $g_i$ can be viewed as central closed even elements of $\Aw$.  Given any dg algebra $A$,  a collection $\{c_j\}_{j=1,\ldots,k}$ of central closed even elements of $A$ gives rise to ``formal deformations''  $\rHN_*(A,\bft)$ and $\rHP_*(A,\bft)$ of $\rHN_*(A)$ and $\rHP_*(A)$, respectively: roughly, these are the negative and the periodic cyclic homology of the formal deformation $(A\zz,z_1c_1+\ldots+z_kc_k)$ of $A$ in the class of {\it curved} dg algebras. The $\C\zz\uu$-module $\rHP_*(A,\bft)$ carries the Gauss-Manin-Getzler connection ${\bconn}_{{\!\!{\partial_{z_i}}}}^{A,\bft}$ (\cite{Get}) and one additional differential operator ${\bconn}_{{\!\!{\partial_{u}}}}^{A,\bft}$, a deformation of ${\bconn}_{{\!\!{\partial_{u}}}}^{A}$. The former operators commute with the latter one; we call ${\bconn}^{A,\bft}$ the  {\it extended} Gauss-Manin-Getzler connection. Furthermore, the pairing $\met_A$ extends to a $\C\zz$-linear $\C\uuu$-sesquilinear $\C\zzu$-valued pairing  $\met_{A, \bft}$ on $\rHN_*(A,\bft)$, and we show that the induced pairing on $\rHP_*(A,\bft)$ is  ${\bconn}^{A,\bft}\,\,$--flat. Finally, we prove that $\bI$ extends to isomorphisms
\[
\bI: \rHN_*(\Aw,\bfg)\simeq \cH^{(0)}_{f,\bfg}, \quad \rHP_*(\Aw,\bfg)\simeq \cH_{f,\bfg}
\]
and that
\[
({\bconn}_{{\!\!{\partial_u}}}^{f,\bfg}-\frac{n}{2u})\cdot\bI=\bI\cdot {\bconn}_{{\!\!{\partial_u}}}^{\Aw,\bfg},\quad {\bconn}_{{\!\!{\partial_{z_i}}}}^{f,\bfg}\cdot\bI=\bI\cdot {\bconn}_{{\!\!{\partial_{z_i}}}}^{\Aw,\bfg}.
\]
These equalities, together with M. Saito's theorem and the fact that  $\met_{\Aw, \bfg}$ is non-zero (as we will explain, the latter follows from the {\it homological smoothness} of $\Aw$), yield Theorem \ref{MC}.

\medskip 
{\bf Conventions.} All our complexes and dg algebras are $\ZZ$-graded and $\K$-linear, unless stated otherwise.  We consider only unital dg algebras. The parity of an element $v$ of a $\ZZ$-graded space is denoted by $|v|$. We will use the bold nabla $\conn$ to denote various ``cochain-level'' connections, such as (\ref{clascon}) and (\ref{exGM}); the induced connections on the cohomology  will be denoted by $\bconn$.

\medskip
{\bf Acknowledgements.} I am grateful to Igor Burban, Bernhard Keller, Alexander Polishchuk, Kyoji Saito, Morihiko Saito, Emanuel Scheidegger, and Christian Sevenheck for discussions and correspondence. I am especially indebted to Morihiko Saito for a careful and detailed response to my query regarding his uniqueness theorem and to Alexander Polishchuk for spotting a mistake in a preliminary version of the paper.  

This research was supported by the ERC Starting Independent Researcher Grant StG No. 204757-TQFT (K.~Wendland PI) and a Fellowship from the Freiburg Institute for Advanced Studies (FRIAS), Freiburg, Germany.

\section{The canonical pairing on the periodic cyclic homology}
\subsection{Flat pairings on mixed complexes} The object that underlies the definition (\ref{tw})  is the  {\it mixed complex} $(\Omega^{\ast}_{\C^n,{\sf 0}}, -d{f}, d)$. The aim of this section is to set up some formalism for dealing with (flat) pairings in the framework of abstract mixed complexes endowed with connections.
 
To begin with, an abstract mixed complex is simply a triple $(\rC_*, b, B)$ where $(\rC_*, b)$ is a complex and $B$ is an odd operator on $\rC_*$ such that $B^2=bB+Bb=0$. (In the above example, the $\ZZ$-grading is given by the parity of the degree of a differential form.)

Let $u$ stand for a formal even variable. Then, obviously, $(\rC_*\uu, b+uB)$ (resp. $(\rC_*\uuu, b+uB)$) is a $\K\uu$-linear (resp. $\K\uuu$-linear) complex. Following \cite{Shk2} we will say that  $(\rC_*, b, B)$ carries a {\it $u$-connection} if  $\rC_*\uu$ is endowed with an (even) operator
\[
{\conn}_{{\!\!{\partial_u}}}=\partial_u+\cA(u),\quad  \cA(u)\in{\End}_{\K\uu}\rC_*\uu
\]
such that
$
[{\conn}_{{\!\!{\partial_u}}},b+uB]=\frac{1}{2u}(b+uB). 
$
This condition implies that ${\conn}_{{\!\!{\partial_u}}}$ induces a connection on $\rH^*(\rC_*\uu, b+uB)$  (which we have agreed to denote by ${\bconn}$).

Furthermore, a {\it morphism} from $(\rC_*, b, B, {\conn}_{{\!\!{\partial_u}}})$ to $(\rC_*', b', B', {\conn}_{{\!\!{\partial_u}}}')$ is a morphism $f(u)$ of the associated $\K\uu$-linear complexes such that
\begin{eqnarray}\label{nablaf}
{\conn}_{{\!\!{\partial_u}}}'f(u)-f(u){\conn}_{{\!\!{\partial_u}}}: (\rC_*\uu, b+uB)\rightarrow (\rC_*'\uu, b'+uB')
\end{eqnarray}
is 0-homotopic. 

\begin{remark}{\rm 
Note that by \cite[Lemma 2.1]{Shk2} the commutator (\ref{nablaf}) is always a morphism of complexes which clarifies the definition.  Also note that in our setting to be 0-homotopic is the same as to induce the trivial map on cohomology since we consider complexes over a field (namely, $\K\uu$).}
\end{remark}

The category of mixed complexes with $u$-connections is monoidal: 
\[
(\rC_*, b, B, {\conn}_{{\!\!{\partial_u}}})\otimes (\rC_*', b', B', {\conn}_{{\!\!{\partial_u}}}'):=(\rC_*\otimes\rC_*', b\otimes1+1\otimes b', B\otimes1+1\otimes B', {\conn}_{{\!\!{\partial_u}}}\otimes1+1\otimes {\conn}_{{\!\!{\partial_u}}}'),
\]
the unit object being $\left(\K, 0, 0, \partial_u\right)$. (Here $\K$ is viewed as a graded space whose even component is $\K$ and the odd component is trivial.)

The last ingredient we need is the following involution on the category of mixed complexes with $u$-connections: 
\[
\left(\rC_*, b, B, {\conn}_{{\!\!{\partial_u}}}\right)^{\vee}:=\left(\rC_*, b, -B, {\conn}_{{\!\!{\partial_u}}}^{\vee}\right)
\]
where $(\partial_u+\cA(u))^\vee:=\partial_u-\cA(-u)$.

Now we are in a position to define pairings on mixed complexes with $u$-connections: By a {\it pairing} on $(\rC_*, b, B, {\conn}_{{\!\!{\partial_u}}})$ we will understand a morphism 
\[
\langle\,\cdot\,,\,\cdot\,\rangle: (\rC_*, b, B, {\conn}_{{\!\!{\partial_u}}})\otimes(\rC_*, b, B, {\conn}_{{\!\!{\partial_u}}})^{\vee}\rightarrow\left(\K, 0, 0, \partial_u\right)
\]

\begin{example}{\rm Let $X$ be a compact complex manifold. Consider the mixed complex $(\mathscr{A}(X),\bar{\partial},{\partial})$ where $\mathscr{A}(X)$ is the $\ZZ$-graded space of complex $C^\infty$-forms on $X$, with the $\ZZ$-grading given by the parity of the degree of differential forms. The mixed complex carries a $u$-connection given by
\[ 
{\conn}_{{\!\!{\partial_u}}}^X=\partial_u+\frac{\gamma'}{u},  \quad \gamma'|_{\mathscr{A}^{p,q}(X)}:=\frac{q-p}2\cdot \mathrm{id}_{\mathscr{A}^{p,q}(X)}
\]
Define 
\[
\langle\omega_1,\omega_2\rangle_X:=\int_X \omega_1\wedge\omega_2^{\vee}
\]
where $\omega^{\vee}:=(-1)^p\cdot\omega$ for $\omega\in\mathscr{A}^{p,\ast}(X)$ (cf. \cite{Ra}). 
It is easy to check that  $\langle\,\cdot\,,\,\cdot\,\rangle_X$ extends to a pairing on $(\mathscr{A}(X),\bar{\partial},{\partial}, {\conn}_{{\!\!{\partial_u}}}^X)$.}
\end{example}

Recall the automorphism $\lambda(u)^\star:=\lambda(-u)$ of $\K\uu$. Given a mixed complex $(\rC_*, b, B)$, the automorphism induces  canonical $\C\uu$-antilinear isomorphisms
\[ (\rC_*\uu, b+uB)\to(\rC_*\uu, b-uB),\quad \rH^\ast(\rC_*\uu, b+uB)\to\rH^\ast(\rC_*\uu, b-uB)\]
which we will also denote by $\,^\star$. Suppose the mixed complex is equipped with a $u$-connection and a  pairing. Consider the $\C\uu$-sesquilinear pairing 
\begin{eqnarray*}
\rC_*\uu\times\rC_*\uu\rightarrow \K\uu, \quad (\rc_1(u), \rc_2(u))\mapsto\langle \rc_1(u), \rc_2(u)^\star\rangle
\end{eqnarray*}
Obviously, it respects the differentials 
\begin{eqnarray*}
\langle (b+uB)\rc_1(u), \rc_2(u)^\star\rangle =-(-1)^{|\rc_1(u)|}\langle \rc_1(u), \left((b+uB)\rc_2(u)\right)^\star\rangle
\end{eqnarray*}
and, hence, descends to $\rH^\ast(\rC_*\uu, b+uB)$. We will denote the induced pairing on the cohomology by $\met$. It follows from the definitions that
\begin{eqnarray*}
\met({\bconn}_{{\!\!{\partial_u}}}\,\,\cdot\,,\, \cdot\,)-\met(\,\cdot\,, {\bconn}_{{\!\!{\partial_u}}}\,\cdot\,)=\partial_u\met(\,\cdot\,,\, \cdot\,).
\end{eqnarray*}

\begin{remark}{\rm  In the next sections, we will have to deal with mixed complexes over the algebra $\C\zz$. Obviously, all the definitions in this section make sense in this setting, and we will be freely using the terminology introduced above.}
\end{remark}

\subsection{The cyclic mixed complex of a dg algebra and the canonical $u$-connection}
In this section we recall the definition and some properties of the canonical connection ${\bconn}^A$ on the periodic cyclic homology $\rHP_*(A)$ of a dg algebra $A$. More details can be found in \cite{Shk2,Shk3}.

Let $A=(A,d)$ be a dg algebra and $(\rCn_*(A), b, B)$ be its (normalized) cyclic mixed complex. Recall that
$
\rCn_*(A):=\oplus_{l\geq0} A\otimes (\Pi\,\overline{A})^{\otimes l}
$
where $\overline{A}:=A/\C$ and $\Pi$ is the parity reversing functor. (We will write elements of 
$
A\otimes (\Pi\,\overline{A})^{\otimes l}
$ 
as $a_0[a_1|\ldots |a_l]$.) The operator $b$ is the Hochschild differential. It is the sum of two anti-commuting differentials $b(\delta)$ and $b(\mu)$ where the former is the total differential on the tensor products of the complexes and the latter is a version of the usual Hochschild differential \cite{Lod} which takes account of the non-trivial $\ZZ$-grading. Specifically, if we set
\begin{eqnarray*}\label{lm} 
&&\ld^{(0)}(a_0[a_1|a_2|\ldots |a_l]):=da_0[a_1|a_2|\ldots|a_l],\quad \lm^{(0)}(a_0[a_1|a_2|\ldots |a_l]):=(-1)^{|a_0|}a_0a_1[a_2|\ldots|a_l],\\
&&\tau(a_0[a_1|a_2|\ldots |a_l]):=(-1)^{(|a_0|+1)(l+\sum_{i=1}^l|a_i|)}a_1[a_2|\ldots |a_l|a_0],\\
&&\lm^{(i)}:=\tau^{-i}\lm^{(0)}\tau^{i},\quad \ld^{(i)}:=\tau^{-i}\ld^{(0)}\tau^{i}
\end{eqnarray*}
then $\bd=\sum_i \ld^{(i)}$ and $\bm=\sum_i \lm^{(i)}$. The operator $B$ is defined as the composition $s\cdot N$ where
\begin{equation}\label{sN}
s(a_0[a_1|a_2|\ldots |a_l]):=1[a_0|a_1|a_2|\ldots |a_l], \quad N|_{A\otimes (\Pi\,\overline{A})^{\otimes l}}:=\sum_{i=0}^l\tau^i.
\end{equation}

The negative (resp. periodic) cyclic homology of $A$ are the cohomology of the associated $\K\uuu$-linear (resp. $\K\uu$-linear) complex:
\[
 \rHN_*(A):=\rH^*(\rCn_* (A)\uuu, b+ uB),\quad \rHP_*(A):=\rH^*(\rCn_* (A)\uu, b+ uB).
\]

Consider the endomorphisms $\gamma$, $\U$, and $\V$ of $\rCn_*(A)$ whose restrictions to ${A\otimes (\Pi\,\overline{A})^{\otimes l}}$ are defined by the formulas
\begin{equation}\label{end1}
\gamma= l\cdot {\rm id}, \quad \U=-\lm^{(0)}\ld^{(1)}, \quad \V=-\sum_{i=1}^{l}\sum_{j=i+1}^{l+1}s\tau^{j}\ld^{(i)}.
\end{equation}
The canonical $u$-connection on $(\rCn_*(A), b, B)$ is given by
\begin{eqnarray*}
{\conn}_{{\!\!{\partial_u}}}^A=\partial_u+\frac{\U}{2u^2}+\frac{\V-\gamma}{2u}
\end{eqnarray*}
Let us recall some of its basics properties established in \cite{Shk2,Shk3}. 

Let $A^{\rm op}$ denote the  opposite dg algebra, i.e.  the same complex $(A,d)$ endowed with the opposite product
$
a_1\otimes a_2\mapsto(-1)^{|a_1||a_2|}a_2a_1.
$
Consider the isomorphism of graded spaces $\Phi: \rCn_*(A)\rightarrow\rCn_*(A^{\rm op})$ given by
\begin{eqnarray}\label{Phi}
\Phi(a_0[a_1|a_2|\ldots |a_l])=(-1)^{l+\sum_{1\leq i<j\leq
l}(|a_i|+1)(|a_j|+1)}a_0[a_l|a_{l-1}|\ldots |a_1]
\end{eqnarray}
The following claim is a consequence of Propositions 3.2, 3.3, and 3.5 in \cite{Shk2}.
\begin{proposition}\label{opp} $\Phi$ induces a morphism from $(\rCn_*(A), b, B, {\conn}_{{\!\!{\partial_u}}}^A)^{\vee}$ to $(\rCn_*(A^{\rm op}), b, B, {\conn}_{{\!\!{\partial_u}}}^{A^{\rm op}})$.
\end{proposition}

\medskip

Let $A'$ be another dg algebra. Consider the map 
$\ck: \rCn_*(A)\otimes \rCn_*(A')\rightarrow\rCn_*(A\otimes A')$
defined by
\begin{multline}\label{ck}
\ck(a_0[a_1|\ldots |a_l]\otimes
a'_0[a'_1|\ldots|a'_m])\\=(-1)^{|a'_0|(l+\sum_{i=1}^l|a_i|)}(a_0\otimes a'_0)[\sh\{a_1\otimes1|\ldots|a_l\otimes1\}\{1\otimes a'_1|\ldots|1\otimes a'_m\}]
\end{multline}
where $\sh\{\ldots\}\{\ldots\}$ stands for the sum over all the permutations that shuffle the $a$-terms with the $a'$-terms while preserving the order of the former and the latter (every transposition
$[\,\ldots|x|y|\ldots\,]\rightarrow[\,\ldots|y|x|\ldots\,]$ contributes $(-1)^{(|x|+1)(|y|+1)}$ to the sign in front of the new tensor).
Furthermore, define  
$\cK: \rCn_*(A)\otimes \rCn_*(A')\rightarrow\rCn_*(A\otimes A')$ by
\begin{multline}\label{cK}
\cK(a_0[a_1|\ldots |a_l]\otimes
a'_0[a'_1|\ldots|a'_m])\\=(-1)^{l+\sum_{i=0}^l|a_i|}(1\otimes 1)[\Sh\{a_0\otimes1|\ldots|a_l\otimes1\}\{1\otimes a'_0|\ldots|1\otimes a'_m\}]
\end{multline}
where $\Sh\{\ldots\}\{\ldots\}$ stands for the sum (again, with signs defined as before) over all the permutations that cyclically permute $a_0, \ldots, a_l$ and $a'_0, \ldots,a'_m$, and then shuffle the $a$-terms with the $a'$-terms so that $a_0$ stays to the left of $a'_0$. As it is explained in \cite{Lod}, the linear combination $\ck+u\cK$ commutes with the cyclic differentials
\begin{equation*} (b+uB)\,(\ck+u\cK)=(\ck+u\cK)\,((b+uB)\otimes 1+1\otimes (b+uB))
\end{equation*} 
and induces K\"unneth-type maps in various versions of cyclic homology. The following fact is  the main result of \cite{Shk3}:

\begin{proposition}\label{kun} The operator $\ck+u\cK$ defines a morphism
\[
(\rCn_*(A), b, B, {\conn}_{{\!\!{\partial_u}}}^A)\otimes (\rCn_*(A'), b, B, {\conn}_{{\!\!{\partial_u}}}^{A'})\rightarrow (\rCn_*(A\otimes A'), b, B, {\conn}_{{\!\!{\partial_u}}}^{A\otimes A'})
\]
\end{proposition}

\subsection{The canonical pairing}\label{11} 

From now on, $A=(A,d)$ stands for a {\it proper} dg algebra, i.e. we assume that the cohomology $\rH A=\rH^*(A,d)$ is finite-dimensional.\footnote{Let us also assume that $\rH A\neq\{0\}$.}
The aim of this section is to describe a pairing on $(\rCn_*(A), b, B, {\conn}_{{\!\!{\partial_u}}}^A)$ that gives rise to the canonical pairing $\met_A$ on $\rHP_*(A)$ defined in (\ref{canpr}). 

As complexes, $A$ and ${\rH}A$ are homotopy equivalent. In what follows, we use some explicit deformation retract data
\begin{eqnarray}\label{drd1}
p: A\rightarrow {\rH}A, \quad i: {\rH}A\rightarrow A,\quad h: A\rightarrow \Pi A\\
pd=0,\quad di=0, \quad pi=\mathrm{id}_{{\rH}A},\quad ip=\mathrm{id}_{A}-[d,h] \nonumber
\end{eqnarray}
which we may and will assume to satisfy the so-called side conditions:
\begin{eqnarray}\label{sc}
ph=0,\quad hi=0,\quad h^2=0
\end{eqnarray}

Let  ${\End}\, A$ be the dg algebra of endomorphism of $A$ (viewed as a complex) and consider the following sequence of linear maps:
\begin{equation*}
\rCn_*(A)\otimes\rCn_*(A)\stackrel{1\otimes\Phi}{\longrightarrow}\rCn_*(A)\otimes\rCn_*(A^{\rm op})\stackrel{\ck}{\longrightarrow}\rCn_*(A\otimes A^{\rm op})\stackrel{\rCn_*(\rho)}{\longrightarrow}\rCn_*({\End}\, A)\stackrel{\STR}{\longrightarrow}\K.
\end{equation*}
Here $\Phi$ and $\ck$ are as in (\ref{Phi}) and (\ref{ck}), $\rCn_*(\rho)$ is the morphism of complexes induced by the homomorphism of dg algebras  
\begin{eqnarray*}
\rho: A\otimes A^{\rm op}\,{\rightarrow}\,{\End}\, A,\quad \rho: a_1\otimes a_2\mapsto L_{a_1}R_{a_2}\\
L_{a_1}: a\mapsto a_1a,\quad R_{a_2}: a\mapsto (-1)^{|a||a_2|}aa_2
\end{eqnarray*}
and $\STR$ is the Feigin-Losev-Shoikhet super-trace defined by the formula (cf. \cite{FLS})
\begin{eqnarray}\label{fls}
\quad\STR(T_0[T_1|\ldots|T_l])=\str_{{\rH}A}\left(F \left(N\left(T_0[T_1|\ldots|T_l]\right)\right)\right)
\end{eqnarray}
where $N$ is the same as in (\ref{sN}), 
\[
F: \rCn_*({\End}\, A)\to {\End}\, {\rH}A,\quad  T_0[T_1|\ldots|T_l]\mapsto pT_0hT_1h\ldots hT_li,
\] and $\str_{{\rH}A}$ is the ordinary super-trace of a linear operator on ${{\rH}A}$. Note that (\ref{fls}) is a well-defined functional on $\rCn_*({\End}\, A)$ thanks to the side conditions (\ref{sc}).

\begin{proposition}\label{thm1}
The composition 
\begin{equation}\label{canpair}
\langle \,\cdot\,,\,\cdot\, \rangle_A:=\STR\cdot\rCn_*(\rho)\cdot\ck\cdot(1\otimes\Phi),
\end{equation}  
when extended by $\K\uu$-linearity, defines a pairing on $(\rCn_*(A), b, B, {\conn}_{{\!\!{\partial_u}}}^A)$.
\end{proposition}
\noindent{\bf Proof.} Observe that  the right-hand side of (\ref{canpair}) equals
$
\STR\cdot\rCn_*(\rho)\cdot(\ck+u\cK)\cdot(1\otimes\Phi)
$
since  the image of the linear map $\rCn_*(\rho)\cdot\cK$ is spanned by tensors of the form 
$
\mathrm{id}_{A}[T_1|\ldots|T_l]\in\rCn_*({\End}\, A)
$ 
(see (\ref{cK})), and these are annihilated by $\STR$ due to  (\ref{sc}).
Thus, in view of Propositions \ref{opp} and \ref{kun}, it suffices to show that $\STR$ induces a morphism \[\left(\rCn_*({\End}\, A), b, B, {\conn}_{{\!\!{\partial_u}}}^{{\End}\, A}\right)\rightarrow\left(\K, 0, 0, \partial_u\right).\]
That $\STR$ is a morphism of the associated $\K\uu$-linear complexes  (i.e. $\STR\cdot(b+uB)=0$) was established in \cite{FLS}. What we need to prove is that the operator
\begin{equation}\label{22}
\partial_u\cdot\STR-\STR\cdot{\conn}_{{\!\!{\partial_u}}}^{{\End}\, A}=-\frac1{2u^2}\STR\cdot (\U+(\V-\gamma)u),
\end{equation}
induces the trivial functional on $\rHP_*({\End}\, A)$. 

In fact,  $\rHP_*({\End}\, A)$ is very simple: by the dg Morita invariance of the cyclic homology \cite{K2,K3}  it is isomorphic to $\rHP_*(\K)$ (that is, it sits in degree 0 and is one-dimensional). Moreover, we can write out an explicit generator:

\begin{lemma}  Let $\pi_0\in{\End}\, \rH A$ be any even rank 1 projector and set ${\pi}:=i\cdot\pi_0\cdot p\in {\End}\, A$.  Then
\begin{equation}\label{gener}
{\pi}+\sum_{l=1}^\infty(-1)^{l}\frac{(2l)!}{l!}u^l({\pi}[\underbrace{{\pi}|\ldots|{\pi}}_{2l\,\,\text{\rm copies}}]-\frac12\cdot\mathrm{id}_{A}[\underbrace{{\pi}|\ldots|{\pi}}_{2l\,\,\text{\rm copies}}])
\end{equation}
is a $(b+uB)$-cocycle and its class in $\rHP_0({\End}\, A)$ is non-trivial.
\end{lemma}
\noindent{\bf Proof of the Lemma.} The element (\ref{gener}) is nothing but the non-commutative Chern character of the idempotent ${\pi}$; cf. \cite{Lod}. Alternatively, one can prove that it is a cocycle directly and then apply $\STR$ to show that it is non-trivial.\hfill $\blacksquare$

\medskip

To finish the proof of Proposition \ref{thm1}, it remains to notice that the generator (\ref{gener}) is annihilated by the functional (\ref{22}). 
\hfill $\blacksquare$

\begin{remark}{\rm
Notice that the isomorphism $\rHP_*({\End}\, A)\simeq\K\uu$ induced by $\STR$ is independent of the deformation retract data. Indeed, for any choice of (\ref{drd1}) $\STR$ sends the element (\ref{gener}) to $\pm1$ where the sign depends on the parity of the image of $\pi$.  As a consequence, the resulting sesquilinear pairing  
$
\rHP_*(A)\times\rHP_*(A)\rightarrow\K\uu
$ 
is also independent of the deformation retract data. It is easy to see this latter pairing is nothing but $\met_A$.
}
\end{remark}

Recall \cite{KS} that a dg algebra $A$ is called {\it homologically smooth} if it is a perfect bimodule over itself. It has been conjectured in \cite{KS} that any such dg algebra possesses the following property known as the non-commutative Hodge--to--de Rham degeneration condition: The negative cyclic homology  $\rHN_*(A):=\rH^*(\rCn_* (A)\uuu, b+ uB)$ of $A$ is a free $\K\uuu$-module; equivalently, the canonical $\K\uuu$-linear map $\rHN_*(A)\to\rHP_*(A)$ is an embedding. Since the conjecture is still open, this condition is added as an extra assumption in the statement of our next result.

\begin{proposition}\label{thm2}
If $A$ is a proper and homologically smooth dg algebra satisfying the Hodge--to--de Rham degeneration condition then $\met_A$ is non-degenerate.
\end{proposition}
\noindent{\bf Sketch of the proof.} The Hochschild homology of a proper homologically smooth dg algebra is finite-dimensional (cf., for instance, \cite{Shk0}). It follows from the Hodge-to-de Rham degeneration condition that $\rHN_*(A)/u\rHN_*(A)\simeq \rHH_*(A)$. Thus,  $\rHN_*(A)$ is a free $\K\uuu$-module whose rank equals the dimension of $\rHH_*(A)$. 

Let us identify $\rHN_*(A)$ with its image in $\rHP_*(A)$ under the above canonical map. 
By (\ref{canpair})  the restriction of $\met_A$ onto $\rHN_*(A)$ takes values in $\K\uuu$. Consider the induced $\K$-valued pairing $\eta_A$ on $\rHH_*(A)$. Obviously, to prove the proposition it is enough to show that $\eta_A$ is non-degenerate. This was proved in the $\Z$-graded case  in \cite[Theorem 5.3]{Shk1} but the argument works in the $\ZZ$-graded setting as well (cf. \cite{PV}). 
\hfill $\blacksquare$

\begin{remark}{\rm Note the higher residue pairing is (skew)-hermitian: 
$
\Kw(\omega_1, \omega_2)=(-1)^n\cdot \Kw(\omega_2, \omega_1)^\star.
$
On the other hand, $\met_A$ is neither hermitian nor skew-hermitian in general. This can be easily seen in the example of the path algebra of any non-trivial finite quiver without oriented cycles. An immediate idea is to fix this by replacing $\met_A$ with its (skew-)hermitization (the resulting pairing is still flat).
Note, however, that the new pairing may not be non-degenerate even if the original one was. 
We believe that $\met_A$ is automatically (skew-)hermitian when $A$ is a {\it Calabi-Yau} dg algebra \cite{KS}. }
\end{remark}

\section{Deformations}
\subsection{Deformation of the canonical $u$-connections}
Given an even element $c\in A$ we will denote by $b(c)$ the odd operator on $\rCn_*(A)$ defined by 
$
b(c)|_{A\otimes (\Pi\,\overline{A})^{\otimes l}}:=\sum_{i=1}^{l+1}c^{(i)}
$
where
\[
c^{(i)}(a_0[a_1|a_2|\ldots |a_l]):=(-1)^{i+\sum_{k=0}^{i-1}|a_k|}a_0[a_1|\ldots|a_{i-1}|c|a_{i}| \ldots|a_l].
\]

Let $\bft=\{c_j\}_{j=1,\ldots,k}$ be a collection of {\it central closed even} elements of $A$. Such a collection gives rise to a deformation of the mixed complex $(\rCn_*(A), b, B)$, namely 

\begin{proposition}\label{le1} Fix $k$ variables ${z}_1,\ldots,{z}_k$. Then
\begin{equation}\label{dmc}
(\rCn_*(A)[\![{z}_1,\ldots,{z}_k]\!],\,\, b+{z}_1b(c_1)+\ldots+{z}_kb(c_k), \,\,B)
\end{equation}
is a $\C\zz$-linear mixed complex.
\end{proposition}
\noindent{\bf Proof.} See, for instance, \cite[Lemma D.1]{Shk2}. \hfill $\blacksquare$

\begin{remark}{\rm This mixed complex may be thought of as the cyclic mixed complex of the $\C\zz$-linear {\it curved} dg algebra  $(A\zz,z_1c_1+\ldots+z_kc_k)$. Note, however, that  \[\rCn_*(A\zz)\varsubsetneq\rCn_*(A)\zz\varsubsetneq\rCp(A\zz).\]}
\end{remark}

Let us introduce notation for the cohomology groups associated with the above mixed complex:
\begin{eqnarray*}
&&\rHN_*(A,\bft):=\rH^*(\rCn_*(A)\zzu,b+\sum{z}_jb(c_j)+uB),\\
&&\rHP_*(A,\bft):=\rH^*(\rCn_*(A)\zz\uu,b+\sum{z}_jb(c_j)+uB).
\end{eqnarray*} 

Associated with an even element $c\in A$ are two more (this time even) operators on $\rCn_*(A)$ whose restrictions to ${A\otimes (\Pi\,\overline{A})^{\otimes l}}$ are defined by the formulas (cf. (\ref{end1}))
\begin{equation*}
e(c)=-\lm^{(0)}c^{(1)}, \quad E(c)=-\sum_{i=1}^{l+1}\sum_{j=i+1}^{l+2}s\tau^{j}c^{(i)}.
\end{equation*}

\begin{proposition} Given a collection $\bft=\{c_j\}$ as above
\begin{equation*}
{\conn}_{{\!\!{\partial_u}}}^{A,\bft}:={\conn}_{{\!\!{\partial_u}}}^A+\frac{{z}_1e(c_1)+\ldots+{z}_ke(c_k)}{u^2}+\frac{{z}_1E(c_1)+\ldots+{z}_kE(c_k)}{u}
\end{equation*}
is a $u$-connection on (\ref{dmc}).
\end{proposition}
\noindent{\bf Proof.} We need to show that
\begin{equation*}
\left[\partial_u+\frac{\U}{2u^2}+\frac{\V-\gamma}{2u}+\frac{\sum z_je(c_j)}{u^2}+\frac{\sum z_jE(c_j)}{u}, b+\sum z_jb(c_j)+uB\right]=\frac{1}{2u}(b+\sum z_jb(c_j)+uB).
\end{equation*}
This equality follows from 
\begin{equation}\label{i}
[{\U}+u{\V}, b+uB]=u b(\delta), \quad [{e(c)}+u{E(c)}, b+uB]=u b(c)
\end{equation}
and
\begin{equation}\label{ii}
[{\U}+u{\V}, b(c)]=0, \quad [{e(c)}+u{E(c)}, b(c')]=0,
\end{equation}
where $c$ and $c'$ are arbitrary central closed even elements. Both formulas in (\ref{i}) are special cases of \footnote{In fact, the former equality is equivalent to the claim that ${\conn}_{{\!\!{\partial_u}}}^A$ is a $u$-connection.} \cite[Eq.(2.1)]{Get}, while (\ref{ii}) can be 
verified using the formulas
\begin{eqnarray*}
\ld^{(i)}c^{(j)}
=\begin{cases}
-c^{(j)}\ld^{(i-1)}\quad (i>j)\\
0\qquad \qquad \quad (i=j)\\
-c^{(j)}\ld^{(i)}\qquad (i<j)\\
\end{cases},
\quad
c^{(i)}{c'}^{(j)}=-
\begin{cases}
{c'}^{(j)}c^{(i-1)} \quad (i>j)\\
{c'}^{(j+1)}c^{(i)} \quad (i\leq j)\\
\end{cases}\\
c^{(i)}\lm^{(0)}=-\lm^{(0)}c^{(i+1)}\quad(i\geq1), \quad c^{(i)}=\tau^{-i}c^{(0)}\tau^{i},\quad  sc^{(i)}=-c^{(i+1)}s.\qquad \blacksquare
\end{eqnarray*}

\subsection{Extended Gauss-Manin-Getzler connection}
Let ${\conn}_{{\!\partial_{z_i}}}^{A,\bft}$ be the operator on $\rCn_*(A)[\![{z}_1,\ldots,{z}_k]\!]\uu$ given by
\begin{equation*}
{\conn}_{{\!\!{\partial_{z_i}}}}^{A,\bft}=\partial_{z_i}-\frac{e(c_i)}{u}-E(c_i).
\end{equation*}

\begin{proposition} The operator ${\conn}_{{\!\!{\partial_{z_i}}}}^{A,\bft}$ descend to
$
\rHP_*(A,\bft).
$
\end{proposition}
\noindent{\bf Proof.}  
\[
[{\conn}_{{\!\!{\partial_{z_i}}}}^{A,\bft}, b+\sum_j{z}_jb(c_j)+uB]=b(c_i)-\frac{1}{u}[{e(c_i)}+u{E(c_i)}, b+uB]-\frac{1}{u}\sum_j{z}_j[{e(c_i)}+u{E(c_i)}, b(c_j)]
\]
which equals 0 by  (\ref{i}) and (\ref{ii}).
\hfill $\blacksquare$

\medskip

Thus, by the two previous propositions, $\rHP_*(A,\bft)$ carries a connection ${\bconn}^{A,\bft}$.

\begin{proposition} The connection ${\bconn}^{A,\bft}$ is flat.
\end{proposition}
\noindent{\bf Proof.} We will show that the commutators $[{\conn}_{{\!\!{\partial_{z_i}}}}^{A,\bft}, {\conn}_{{\!\!{\partial_{z_j}}}}^{A,\bft}]$ and  $[{\conn}_{{\!\!{\partial_{u}}}}^{A,\bft},{\conn}_{{\!\!{\partial_{z_i}}}}^{A,\bft}]$ are homotopic to 0. The claim is a special case\footnote{Formally speaking, Theorem 3.3 of \cite{Get} applies to $[{\conn}_{{\!\!{\partial_{z_i}}}}^{A,\bft}, {\conn}_{{\!\!{\partial_{z_j}}}}^{A,\bft}]$ only since the operator ${\conn}_{{\!\!{\partial_u}}}^{A,\bft}$ was not considered in \cite{Get}. However, the formulas needed to prove the claim for $[{\conn}_{{\!\!{\partial_{u}}}}^{A,\bft},{\conn}_{{\!\!{\partial_{z_i}}}}^{A,\bft}]$ can be found in {\it loc. cit.} as well.} of \cite[Theorem 3.3]{Get}, so we will only sketch the proof.
 
To prove that $[{\conn}_{{\!\!{\partial_{z_i}}}}^{A,\bft}, {\conn}_{{\!\!{\partial_{z_j}}}}^{A,\bft}]$ is 0-homotopic, it suffices to show that 
$
[{e(c)}+{u}E(c), {e(c')}+{u}E(c')]
$
is 0-homotopic for any pair of central closed even elements $c,c'\in A$. By definition,
\[
e(c)(a_0[a_1|\ldots |a_l])=a_0c[a_1|\ldots|a_l],
\]
which implies $[{e(c)}, {e(c')}]=0$. Also, $E(c)E(c')=E(c')E(c)=0$ by the definition of the operators (recall that we work with the normalized Hochschild complex). Thus, the only non-trivial part of the claim is that
$
[{e(c)},E(c')]+[E(c),{e(c')}]
$
is 0-homotopic. One can show that 
\begin{multline*}
[{e(c)},E(c')]+[E(c),{e(c')}]+b(c)b(c') =(b+\sum z_jb(c_j)+uB) H(c,c')+H(c,c')(b+\sum z_jb(c_j)+uB) 
\end{multline*}
where $H(c,c')$ is an odd operator on $\rCn_*(A)$ whose restriction to ${A\otimes (\Pi\,\overline{A})^{\otimes l}}$ is given by 
\[
H(c,c')=\sum_{i=1}^{l+1}\sum_{j=i+1}^{l+2}\sum_{m=j+1}^{l+3}s\tau^m(c^{(j)}c'^{(i)}-c'^{(j)}c^{(i)}).
\]
It remains to observe that the formulas (\ref{i}) and (\ref{ii}) imply that  $b(c)b(c')$ is 0-homotopic.

Let us explain why $[{\conn}_{{\!\!{\partial_{u}}}}^{A,\bft},{\conn}_{{\!\!{\partial_{z_i}}}}^{A,\bft}]$ is 0-homotopic. It is easy to see that 
\[
[{\conn}_{{\!\!{\partial_{u}}}}^{A,\bft},{\conn}_{{\!\!{\partial_{z_i}}}}^{A,\bft}]=-\frac1{2u^3}[{\U}+{u}\V, {e(c_i)}+{u}E(c_i)]-\frac1{u^3}\sum_jz_j[{e(c_j)}+{u}E(c_j), {e(c_i)}+{u}E(c_i)].
\]
Thus, in view of the previous discussion, it remains to prove that
$
[{\U}+{u}\V, {e(c)}+{u}E(c)]
$
is 0-homotopic for any central closed even element $c\in A$. The only non-trivial part of the claim is that
$
[{\U},E(c)]+[\V,{e(c)}]
$
is 0-homotopic. One shows that 
\begin{eqnarray*}
[{\U},E(c)]+[\V,{e(c)}]+\bd b(c) =(b+\sum z_jb(c_j)+uB) H(c)+H(c)(b+\sum z_jb(c_j)+uB) 
\end{eqnarray*}
where $H(c)$ is an odd operator on $\rCn_*(A)$ whose restriction to ${A\otimes (\Pi\,\overline{A})^{\otimes l}}$ is given by 
\[
H(c)=\sum_{i=1}^{l}\sum_{j=i+1}^{l+1}\sum_{m=j+1}^{l+2}s\tau^m(\delta^{(j)}c^{(i)}-c^{(j)}\delta^{(i)}).
\]
Again,  it follows from (\ref{i}) and (\ref{ii})  that  $\bd b(c)$ is 0-homotopic.
\hfill $\blacksquare$

\subsection{Extension of the canonical pairing} Let $A$ be a proper dg algebra. Let 
\[
\langle \,\cdot\,,\,\cdot\, \rangle_{A,\bft}: (\rCn_*(A)\otimes \rCn_*(A))[\![{z}_1,\ldots,{z}_k]\!]\to\C[\![{z}_1,\ldots,{z}_k]\!].
\]
be the $\C\zz$-linear extension of the pairing $\langle \,\cdot\,,\,\cdot\, \rangle_A$.

\begin{proposition}\label{pairext0}
$\langle \,\cdot\,,\,\cdot\, \rangle_{A,\bft}$ is a pairing on 
\begin{equation}\label{mc0}
(\rCn_*(A)[\![{z}_1,\ldots,{z}_k]\!],\,\, b+{z}_1b(c_1)+\ldots+{z}_kb(c_k), \,\,B,\,\, {\conn}_{{\!\!{\partial_u}}}^{A,\bft}).
\end{equation} 
\end{proposition}
\noindent{\bf Proof} is given in Appendix \ref{B}. \hfill$\blacksquare$

\medskip

Let $\met_{A,\bft}$  stand for the corresponding ($\C\zz$-linear $\C\uu$-sesquilinear) pairing on $\rHP_*(A,\bft)$.

\begin{proposition}\label{flat2} The pairing $\met_{A,\bft}$ is ${\bconn}^{A,\bft}\,\,$--flat:
\begin{eqnarray*} 
\met_{A,\bft}({\bconn}_{{\!\!{\partial_{u}}}}^{A,\bft}\,\,\cdot\,,\,\cdot\,)-\met_{A,\bft}(\,\cdot\,,{\bconn}_{{\!\!{\partial_{u}}}}^{A,\bft}\,\,\cdot\,)=\partial_u\met_{A,\bft}(\,\cdot\,,\,\cdot\,),\\
\met_{A,\bft}({\bconn}_{{\!\!{\partial_{z_i}}}}^{A,\bft}\,\,\cdot\,,\,\cdot\,)+\met_{A,\bft}(\,\cdot\,,{\bconn}_{{\!\!{\partial_{z_i}}}}^{A,\bft}\,\cdot\,)=\partial_{z_i}\met_{A,\bft}(\,\cdot\,,\,\cdot\,).
\end{eqnarray*}
\end{proposition}
\noindent{\bf Proof.}
The first equality is a corollary of Proposition \ref{pairext0} while the second  one follows from (\ref{phi}) and Lemma \ref{etEt}.
\hfill$\blacksquare$

\section{Proof of Theorem \ref{MC}}
\subsection{The dg algebra $\Aw$ and the isomorphism $\bI$}\label{conj}
Let us recall the description of the dg algebra $\Aw$ from \cite{Dyck}. We need to fix a decomposition
$f=x_1f_1+\ldots+x_nf_n$ where $f_i$ are some polynomials. Then, as a graded algebra  
$
A_f=\C[\![x_1,\ldots,x_n]\!]\otimes{\End} P_n
$
where $P_n$ is the space of polynomials in $n$ odd variables $\theta_1,\ldots,\theta_n$. The differential of $\Aw$ is defined as the super-commutator with $\Dw:=\sum_i\left(f_i\otimes\theta_i+x_i\otimes {\partial_{\theta_i}}\right)\in\Aw$.

Let us recall now an explicit formula for the isomorphism(s) (\ref{qis}) from \cite{Seg}. More precisely, we will present a formula for a quasi-isomorphism 
\begin{equation}\label{qis1}
\cI: (\rCn_*(\Aw)\uu, b+uB)\to(\widehat{\Omega}^\ast\uu, -d{f}+ud)
\end{equation}
which induces $\bI$ (under the identifications (\ref{ident})).

The map (\ref{qis1}) is defined as the ($\,\K\uu$-linear extension of the) composition 
\begin{equation*}
\rCn_\ast(\C[\![x_1,\ldots,x_n]\!]\otimes{\End} P_k)\xrightarrow{\mathrm{exp}(-b(\Dw))}\rCp(\C[\![x_1,\ldots,x_n]\!]\otimes{\End} P_k)\stackrel{\str}\longrightarrow\rCp(\C[\![x_1,\ldots,x_n]\!])\stackrel{\varepsilon}\longrightarrow\widehat{\Omega}^*
\end{equation*}
where 
\begin{itemize}
\item $\rCp(A):=\prod_{l\geq0} (A\otimes (\Pi\,\overline{A})^{\otimes l})^{\rm even}\bigoplus \prod_{l\geq0} (A\otimes (\Pi\,\overline{A})^{\otimes l})^{\rm odd}$;\item $
b(\Dw)(a_0[a_1|\ldots |a_l]):=\sum_{i=1}^{l+1}a_0[a_1|\ldots|a_{i-1}|\Dw|a_{i}| \ldots|a_l]; 
$
\item 
$
\str (\phi_0\otimes T_0[\phi_1\otimes T_1|\ldots |\phi_l\otimes T_l]):= (-1)^{\sum_{i\,{\rm odd}}|T_i|}\str(T_0T_1\ldots T_l)\phi_0[\phi_1|\ldots |\phi_l] 
$\\ 
where $\phi_i\in \widehat{\Omega}^*$, $T_i\in{\End} P_k$;
\item
 $\varepsilon$ is the Hochschild-Kostant-Rosenberg map: 
\[
\varepsilon(\phi_0[\phi_1|\ldots |\phi_l])=\frac{1}{l!}\phi_0d\phi_1\wedge\ldots \wedge d\phi_l.
\]
\end{itemize}

The main result of \cite{Shk2}  says that $\cI$ defines a morphism of mixed complexes with $u$-connections 
\begin{equation*}
\left(\rCn_*(\Aw), \,\,b, \,\,B, \,\,{\conn}_{{\!\!{\partial_u}}}^{\Aw}\right)\to \left(\widehat{\Omega}^*,\,\, -d{f}, \,\,d, \,\,{\conn}_{{\!\!{\partial_u}}}^{f}-\frac{\gamma}{2u}\right)
\end{equation*}
where ${\conn}_{{\!\!{\partial_u}}}^{f}$ is the operator (\ref{clascon}) and $\gamma|_{\widehat{\Omega}^k}=k\cdot id$.

\subsection{The proof of the theorem}\label{MT}
Let us pick a collection $\bfg=\{g_1,\ldots, g_k\}$ in $\C\xx$ as in the Introduction, i.e. so that the corresponding classes in the Milnor algebra of $f$ form a basis. We will view $\{g_j\}$ as central closed even elements of $\Aw$. Associated with this collection is the $\C\zz$-linear mixed complex\footnote{In this section, $d$ stands for the relative de Rham differential; it commutes with the $z$ variables.}
\begin{eqnarray}\label{gmc}
(\widehat{\Omega}^*\zz,\,\, -d{f}-{z}_1dg_1-\ldots-{z}_kdg_k, \,\,d).
\end{eqnarray}
Let ${\conn}_{{\!\!{\partial_u}}}^{f,\bfg}$ and ${\conn}_{{\!\!{\partial_{z_i}}}}^{f,\bfg}$ denote the operators  (\ref{exGM}). Note that  ${\conn}_{{\!\!{\partial_u}}}^{f,\bfg}-\frac{\gamma}{2u}$ is a $u$-connection on (\ref{gmc}).

As explained in the Introduction, the proof of Theorem \ref{MC} consists in matching all the data
\[
(\widehat{\Omega}^*\zz,\,\, -d{f}-{z}_1dg_1-\ldots-{z}_kdg_k, \,\,d,\,\, {\conn}_{{\!\!{\partial_u}}}^{f,\bfg}-\frac{\gamma}{2u}),\quad \{{\conn}_{{\!\!{\partial_{z_i}}}}^{f,\bfg}\}_{i=1,\ldots,k},
\]
with the similar data associated with $\Aw$
\[
(\rCn_*(\Aw)\zz,\,\, b+{z}_1b(g_1)+\ldots+{z}_kb(g_k), \,\,B,\,\, {\conn}_{{\!\!{\partial_u}}}^{\Aw,\bfg}),\quad \{{\conn}_{{\!\!{\partial_{z_i}}}}^{\Aw,\bfg}\}_{i=1,\ldots,k}.
\]
More precisely, we will show that $\cI$, extended by $\C\zz$-linearity, transforms the latter into the former.

\begin{lemma} The morphism $\cI$ is a morphism between the above mixed complexes; the induced maps 
\[
\bI: \rHN_n(\Aw,\bfg)\to\cH^{(0)}_{f,\bfg},\quad \rHP_n(\Aw,\bfg)\to\cH_{f,\bfg}
\]
 are isomorphisms.
\end{lemma}
\noindent{\bf Proof.} Since we already know that $\cI$ transforms $b$ into $-d{f}$ and $B$ into $d$, we only need to show that for any $g\in\C\xx$ 
\[
\cI\cdot b(g)=-dg\cdot\cI.
\]
Recall that $\cI=\varepsilon\cdot\str\cdot\mathrm{exp}(-b(\Dw))$. Then the claim follows from
\begin{equation}\label{bdbt}
b(\Dw)\cdot b(g)=b(g)\cdot b(\Dw),\quad \str\cdot b(g)=b(g)\cdot \str,\quad \varepsilon\cdot b(g)=-dg\cdot \varepsilon.
\end{equation}
All of these equalities are quite straightforward. Let us, for example, prove the third one:
\begin{eqnarray*}
\varepsilon\cdot b(g)(\phi_0[\phi_1|\ldots |\phi_l])=\sum_{i=1}^{l+1}(-1)^i\varepsilon (\phi_0[\phi_1|\ldots|\phi_{i-1}|g|\ldots|\phi_l])\\
=\frac1{(l+1)!}\sum_{i=1}^{l+1} (-1)^i \phi_0 d\phi_1\wedge\ldots\wedge d\phi_{i-1}\wedge dg\wedge\ldots\wedge d\phi_l\\
=-\frac1{(l+1)!}\sum_{i=1}^{l+1}dg\wedge (\phi_0\wedge d\phi_1\wedge\ldots\wedge d\phi_l)\\
=-dg \wedge (\frac1{l!}\phi_0\wedge d\phi_1\wedge\ldots\wedge d\phi_l)=-dg\wedge\varepsilon(\phi_0[\phi_1|\ldots |\phi_l]).
\end{eqnarray*}

That the extended $\bI$ is still an isomorphism may be shown by induction using, for instance, \cite[Proposition 2.4]{GJ}. 
\hfill$\blacksquare$

\begin{lemma}\label{flat3} The operators
\[
({\conn}_{{\!\!{\partial_u}}}^{f,\bfg}-\frac{\gamma}{2u})\cdot\cI-\cI\cdot {\conn}_{{\!\!{\partial_u}}}^{\Aw,\bfg},\quad {\conn}_{{\!\!{\partial_{z_i}}}}^{f,\bfg}\cdot\cI-\cI\cdot {\conn}_{{\!\!{\partial_{z_i}}}}^{\Aw,\bfg}, 
\]
viewed as morphisms from  
\[
(\rCn_*(\Aw)\zz\uu,\,\, b+\sum{z}_jb(g_j)+uB)\to
(\widehat{\Omega}^*\zz\uu,\,\, -d{f}-\sum{z}_jdg_j+ud),
\]
are 0-homotopic. As a consequence
\[
({\bconn}_{{\!\!{\partial_u}}}^{f,\bfg}-\frac{n}{2u})\cdot\bI=\bI\cdot {\bconn}_{{\!\!{\partial_u}}}^{\Aw,\bfg},\quad {\bconn}_{{\!\!{\partial_{z_i}}}}^{f,\bfg}\cdot\bI=\bI\cdot {\bconn}_{{\!\!{\partial_{z_i}}}}^{\Aw,\bfg}.
\]
\end{lemma}
\noindent{\bf Proof.} It follows from calculations in \cite{Shk2} that
\begin{eqnarray}\label{homotconn}
({\conn}_{{\!\!{\partial_u}}}^{f}-\frac{\gamma}{2u})\cdot\cI - \cI\cdot {\conn}_{{\!\!{\partial_u}}}^{\Aw}
=(-d{f}+ud)\cdot \frac{H}{2u^2}+\frac{H}{2u^2}\cdot(b+uB)
\end{eqnarray}
where
\begin{eqnarray*}
H:=\cI\cdot(\lm^{(0)}+usN')-\varepsilon\cdot\str\cdot(\lm^{(0)}+usN')\cdot\mathrm{exp}(-b(\Dw))- f\cdot d\cdot\cI.
\end{eqnarray*}
where $N'|_{{A\otimes (\Pi\,\overline{A})^{\otimes l}}}:=\sum_{i=1}^{l+1} i\tau^i.$
We will show that
\begin{multline*}
({\conn}_{{\!\!{\partial_u}}}^{f,\bfg}-\frac{\gamma}{2u})\cdot\cI - \cI\cdot {\conn}_{{\!\!{\partial_u}}}^{\Aw,\bfg}=\\
=(-d{f}-\sum z_jdg_j+ud)\cdot \frac{H-\sum z_jg_j\cdot d\cdot\cI}{2u^2}+\frac{H-\sum z_jg_j\cdot d\cdot\cI}{2u^2}\cdot(b+\sum z_jb(g_j)+uB).
\end{multline*}
This equality is a consequence of (\ref{homotconn}), together with the following two formulas: 
\begin{eqnarray}\label{rts1}
-dg\cdot \frac{H}{2u^2}+\frac{H}{2u^2}\cdot b(g)=0
\end{eqnarray}
and 
\begin{eqnarray}\label{rts3}
\,\\ \frac{g}{u^2}\cdot\cI - \cI\cdot \frac{e(g)+uE(g)}{u^2}=
(-d{f}-dg'+ud)\cdot \frac{-g\cdot d\cdot\cI}{2u^2}+\frac{-g\cdot d\cdot\cI}{2u^2}\cdot (b+b(g')+uB)\nonumber
\end{eqnarray}
where $g,g'\in\C[\![{z}_1,\ldots,{z}_k,x_1,\ldots,x_n]\!]$. Let us sketch the proof of these formulas.

One can show that
\begin{eqnarray}\label{bdet}
[b(\Dw), e(g)+uE(g)]=0,\quad [\str, e(g)+uE(g)]=0, \\
\varepsilon\cdot e(g)=g\cdot\varepsilon, \quad [b(g), \mu^{(0)}+usN']=-(e(g)+uE(g)).\nonumber
\end{eqnarray}
These formulas, along with (\ref{bdbt}), yield
\begin{eqnarray*}
-dg\cdot H+H\cdot b(g)
=-dg\cdot\cI\cdot(\lm^{(0)}+usN')+\cI\cdot(\lm^{(0)}+usN')\cdot b(g)\\
+dg\cdot \varepsilon\cdot\str\cdot(\lm^{(0)}+usN')\cdot\mathrm{exp}(-b(\Dw))-\varepsilon\cdot\str\cdot(\lm^{(0)}+usN')\cdot\mathrm{exp}(-b(\Dw))\cdot b(g)\\
+dg\cdot f\cdot d\cdot\cI-f\cdot d\cdot\cI\cdot b(g)\\
=\cI\cdot[b(g),\lm^{(0)}+usN']
-\varepsilon\cdot\str\cdot[b(g),\lm^{(0)}+usN']\cdot\mathrm{exp}(-b(\Dw))
+ f\cdot [dg,d]\cdot\cI\\
=\varepsilon\cdot\str\cdot(\mathrm{exp}(-b(\Dw))\cdot(e(g)+uE(g))-
(e(g)+uE(g))\cdot\mathrm{exp}(-b(\Dw)))=0
\end{eqnarray*}
which proves (\ref{rts1}). 

Furthermore, $dg'\cdot g\cdot d\cdot\cI-g\cdot d\cdot\cI\cdot b(g')=g\cdot [dg',d]\cdot\cI=0$. Therefore, in view of (\ref{bdet}),  equality (\ref{rts3}) would follow from
\begin{eqnarray*}
 2u \cI\cdot E(g)
=(-d{f}+ud)\cdot g\cdot d\cdot\cI+g\cdot d\cdot\cI\cdot (b+uB).
\end{eqnarray*}
The proof of the latter equality is identical to the proof of \cite[Proposition 4.5]{Shk2}.

Note that (\ref{rts3}) also implies that ${\conn}_{{\!\!{\partial_{z_i}}}}^{f,\bfg}\cdot\cI-\cI\cdot {\conn}_{{\!\!{\partial_{z_i}}}}^{\Aw,\bfg}$ is 0-homotopic.
\hfill $\blacksquare$

\medskip

To complete the proof of Theorem \ref{MC}, observe that by Proposition \ref{flat2} and Lemma \ref{flat3} the pairing $u^n\met_{\Aw,\bfg}(\,I_f^{-1}(\cdot)\,,\,I_f^{-1}(\cdot)\,)$ on $\cH_{f,\bfg}$ satisfies the flatness condition (\ref{flat1}). Note also that this pairing is non-trivial. Indeed, $\Aw$ is homologically smooth (see \cite{Dyck}) and therefore, by Proposition \ref{thm2}, the pairing $\met_\Aw$ is non-degenerate which implies $\met_{\Aw,\bfg}\neq0$. Using the uniqueness theorem of M. Saito (see the Introduction) we conclude that
\[
\met_{f,\bfg}=\mathsf{const}\cdot u^n\met_{\Aw,\bfg}(\,I_f^{-1}(\cdot)\,,\,I_f^{-1}(\cdot)\,)
\]
for some constant $\mathsf{const}$. The assertion of Theorem \ref{MC} follows now by `setting' all the $z_i$ to 0.

\appendix

\section{}\label{B}

We will start by proving the following extension of Proposition \ref{opp}:
\begin{lemma}\label{exopp} $\Phi$ induces a morphism (of mixed complexes with $u$-connections) from 
\[\left(\rCn_*(A)\zz,\,\, b+\sum z_jb(c_j), \,\,B,\,\, {\conn}_{{\!\!{\partial_u}}}^{A,\bft}\right)
^{\vee}\]
to
\begin{equation}\label{modcon}
\left(\rCn_*(A^{\rm op})\zz,\,\, b-\sum z_jb(c_j), \,\,B,\,\, {\conn}_{{\!\!{\partial_u}}}^{A^{\rm op}}-\frac{\sum z_je(c_j)}{u^2}-\frac{\sum z_j\widetilde{E}(c_j)}{u}\right)
\end{equation}
where $\widetilde{E}(c):=E(c)+Bb(c)$.
\end{lemma}
\noindent{\bf Proof.} Let us first recall the idea of the proof of Proposition \ref{opp}. Consider the endomorphisms $\widetilde{e}(\delta)$ and $\widetilde{E}(\delta)$ of $\rCn_*(A)$ whose restrictions to ${A\otimes (\Pi\,\overline{A})^{\otimes l}}$ are defined by the formulas
\[
\widetilde{e}(\delta)=\lm^{(l)}\ld^{(l)}, \quad \widetilde{E}(\delta)=\sum_{i=1}^{l}\sum_{j=1}^{i}s\tau^{j}\ld^{(i)}.
\]
One shows that
\[
\Phi^{-1}({\conn}_{{\!\!{\partial_u}}}^{A^{\rm op}})^{\vee} \Phi=\partial_u+\frac{\widetilde{e}(\delta)}{2u^2}+\frac{\widetilde{E}(\delta)-\gamma}{2u}.
\]
and 
\begin{equation*}
2u^2({\conn}_{{\!\!{\partial_u}}}^A-\Phi^{-1}({\conn}_{{\!\!{\partial_u}}}^{A^{\rm op}})^{\vee}\Phi)= (\U-\widetilde{e}(\delta))+u(\V-\widetilde{E}(\delta))=(b+uB)H+H(b+uB)
\end{equation*}
where
\begin{equation}\label{homot1}
H=\ld^{(0)}(1-N).
\end{equation}

To prove  Lemma \ref{exopp}, note that for a central $c$
\begin{eqnarray}\label{phi}
\Phi b(c)=-b(c)\Phi, \quad \Phi e(c)=e(c)\Phi, \quad \Phi E(c)=-\widetilde{E}(c)\Phi
\end{eqnarray}
and therefore
\begin{equation*}
\Phi^{-1}\left({\conn}_{{\!\!{\partial_u}}}^{A^{\rm op}}-\frac{\sum z_je(c_j)}{u^2}-\frac{\sum z_j\widetilde{E}(c_j)}{u}\right)^{\vee}\Phi=\Phi^{-1}({\conn}_{{\!\!{\partial_u}}}^{A^{\rm op}})^{\vee}\Phi+\frac{\sum z_je(c_j)}{u^2}+\frac{\sum z_j{E}(c_j)}{u}.
\end{equation*} 
Thus, the difference between ${\conn}_{{\!\!{\partial_u}}}^{A,\bft}$ and the above $u$-connection is still $(\U-\widetilde{e}(\delta))+u(\V-\widetilde{E}(\delta))$.
One can check that this expression is equal to
$(b+\sum z_jb(c_j)+uB)H+H(b+\sum z_jb(c_j)+uB)$ where $H$ is the same as in (\ref{homot1}) (in other words, $b(c)H+Hb(c)=0$ for any central closed even $c$). 
\hfill $\blacksquare$

\medskip

The above Lemma reduces proving Proposition \ref{pairext0} to showing that the composition $\STR\cdot\rCn_*(\rho)\cdot\ck$ induces a morphism  from the tensor product of
(\ref{mc0}) and (\ref{modcon}) 
to $\left(\C\zz, 0,0,\partial_u\right)$.

\begin{lemma}\label{etEt} For any central closed even $c$
\begin{equation}\label{bt} 
\STR\cdot\rCn_*(\rho)\cdot\ck\cdot(b(c)\otimes1-1\otimes b(c))=0,
\end{equation}
\begin{equation}\label{eE} 
\STR\cdot\rCn_*(\rho)\cdot\ck\cdot(e(c)\otimes1-1\otimes e(c))=0,\quad \STR\cdot\rCn_*(\rho)\cdot\ck\cdot(E(c)\otimes1-1\otimes \widetilde{E}(c))=0.
\end{equation}
\end{lemma}
\noindent{\bf Proof.} Since $\STR\cdot\rCn_*(\rho)\cdot\ck$ is compatible with both $b$ and $B$, (\ref{bt}) is a consequence of  (\ref{eE}) and the second equality in (\ref{i}) (note that  (\ref{i}) remains valid if $E(c)$ is replaced with $\widetilde{E}(c)$). 

The first equality in (\ref{eE}) is quite straightforward (it holds true without $\STR$ in the left-hand side); let us sketch the proof of the second one.

Recall that $\STR=\str_{{\rH}A}\cdot F\cdot N$. We will show that
\[
N\cdot\rCn_*(\rho)\cdot\ck\cdot(E(c)\otimes1)=N\cdot\rCn_*(\rho)\cdot\ck\cdot(1\otimes \widetilde{E}(c)).
\]
Let us fix two tensors
\[
a_0[a_1|\ldots |a_l]\in\rCn_*(A),\quad a'_0[a'_1|\ldots |a'_m]\in\rCn_*(A^{\rm op}).
\]
To simplify the exposition, we will assume that all the $a_j$ and $a'_j$ are {\it odd} elements. (This assumption makes keeping track of signs quite easy.)
We have
\begin{eqnarray}\label{exp1}
N\cdot\rCn_*(\rho)\cdot\ck\cdot(E(c)(a_0[a_1|\ldots |a_l])\otimes
a'_0[a'_1|\ldots|a'_m])\nonumber\\=\sum_{i=1}^{l+1}\sum_{j=i+1}^{l+2} N\left(R_0[\sh\{\tau^{j}(L_0|\ldots|L_{i-1}|\widehat{c}|L_i|\ldots|L_{l})\}\{R_1|\ldots|R_m\}]\right)
\end{eqnarray}
and 
\begin{eqnarray}\label{exp2}
N\cdot\rCn_*(\rho)\cdot\ck\cdot(a_0[a_1|\ldots |a_l]\otimes
\widetilde{E}(c)(a'_0[a'_1|\ldots|a'_m]))\nonumber\\=\sum_{p=1}^{m+1}\sum_{q=1}^{p}  N\left(L_0[\sh\{L_1|\ldots|L_l\}\{\tau^q(R_0|\ldots|R_{p-1}|\widehat{c}|R_p|\ldots|R_m)\}]\right)
\end{eqnarray}
where $L_j:=L_{a_j}$ and $R_j:=R_{a'_j}$, and $\widehat{c}:=L_c=R_c$ (recall that $c$ is central). 
Due to the presence of the operator $N$ in both expressions, it would suffice to compare the parts of the expressions that contain the terms having the same initial element. For instance, let us compare the terms that start with $R_0$. Obviously, in (\ref{exp1}) such terms are
\begin{equation}\label{fe}
\sum_{i=1}^{l+1}\sum_{j=i+1}^{l+2} R_0[\sh\{\tau^{j}(L_0|\ldots|L_{i-1}|\widehat{c}|L_i|\ldots|L_{l})\}\{R_1|\ldots|R_m\}]
\end{equation}
Furthermore, observe that, in general,
\begin{eqnarray}\label{gf1}
\sh\{x_1|\ldots|x_l\}\{y_1|\ldots|y_{k-1}|y|y_k|\ldots|y_m\}=\nonumber\\
\sum_{r=1}^{l+1}\sh\{x_1|\ldots|x_{r-1}\}\{y_1|\ldots|y_{k-1}\}|y|\sh\{x_{r}|\ldots|x_l\}\{y_k|\ldots|y_m\}
\end{eqnarray}
and
\begin{eqnarray}\label{gf2}
\sh\{y_1|\ldots|y_{k-1}|y|y_k|\ldots|y_m\}\{x_1|\ldots|x_l\}=\nonumber\\
\sum_{r=1}^{l+1}\sh\{y_1|\ldots|y_{k-1}\}\{x_1|\ldots|x_{r-1}\}|y|\sh\{y_k|\ldots|y_m\}\{x_{r}|\ldots|x_l\}
\end{eqnarray}
provided all the $x_j$ are odd. Hence,
\begin{eqnarray*}
L_0[\sh\{L_1|\ldots|L_l\}\{\tau^q(R_0|\ldots|R_{p-1}|\widehat{c}|R_p|\ldots|R_m)\}]\\
=L_0[\sh\{L_1|\ldots|L_l\}\{R_q|\ldots|R_{p-1}|\widehat{c}|R_p|\ldots|R_m|R_0|\ldots|R_{q-1}\}]\\
=\sum_{r=1}^{l+1} L_0[\sh\{L_1|\ldots|L_{r-1}\}\{R_q|\ldots|R_{p-1}|\widehat{c}|R_p|\ldots|R_m\}|R_0|\sh\{L_{r}|\ldots|L_l\}\{R_1|\ldots|R_{q-1}\}].
\end{eqnarray*}
Thus, the part of (\ref{exp2}) comprised of the terms with $R_0$ on the first position is 
\begin{eqnarray*}
\sum_{p=1}^{m+1}\sum_{q=1}^{p} \sum_{r=1}^{l+1} R_0[\sh\{L_{r}|\ldots|L_l\}\{R_1|\ldots|R_{q-1}\}|L_0|\sh\{L_1|\ldots|L_{r-1}\}\{R_q|\ldots|R_{p-1}|\widehat{c}|R_p|\ldots|R_m\}].
\end{eqnarray*}
Applying (\ref{gf1}) and (\ref{gf2}) we obtain 
{\small
\begin{eqnarray*}
\sum_{p,q,r} R_0[\sh\{L_{r}|\ldots|L_l\}\{R_1|\ldots|R_{q-1}\}|L_0|\underbrace{\sh\{L_1|\ldots|L_{r-1}\}\{R_q|\ldots|R_{p-1}|\widehat{c}|R_p|\ldots|R_m\}}_{\text{apply (\ref{gf1})}}]=\\
\sum_{p,q,r,s}R_0[\underbrace{\sh\{L_{r}|\ldots|L_l\}\{R_1|\ldots|R_{q-1}\}|L_0|\sh\{L_1|\ldots|L_{s-1}\}\{R_q|\ldots|R_{p-1}\}}_{\text{apply (\ref{gf2})}}|\widehat{c}|\sh\{L_{s}|\ldots|L_{r-1}\}\{R_p|\ldots|R_m\}]\\
=\sum_{p,r,s} R_0[\underbrace{\sh\{L_{r}|\ldots|L_l|L_0|L_1|\ldots|L_{s-1}\}\{R_1|\ldots|R_{p-1}\}|\widehat{c}|\sh\{L_{s}|\ldots|L_{r-1}\}\{R_p|\ldots|R_m\}}_{\text{apply (\ref{gf2})}}]\\
=\sum_{r=1}^{l+1}\sum_{s=1}^{r}  R_0[\sh\{L_{r}|\ldots|L_l|L_0|L_1|\ldots|L_{s-1}|\widehat{c}|L_{s}|\ldots|L_{r-1}\}\{R_1|\ldots|R_m\}]=\text{(\ref{fe})}.\quad\blacksquare
\end{eqnarray*}
}

To complete the proof of Proposition \ref{pairext0}, it remains to prove
\begin{lemma} The map
\begin{eqnarray*}
&&\partial_u\cdot\STR\cdot\rCn_*(\rho)\cdot\ck-\\&&-\STR\cdot\rCn_*(\rho)\cdot\ck\cdot\left( ({\conn}_{{\!\!{\partial_u}}}^{A}+\frac{\sum z_je(c_j)}{u^2}+\frac{\sum z_jE(c_j)}{u})\otimes 1+ 1\otimes({\conn}_{{\!\!{\partial_u}}}^{A^{\rm op}}-\frac{\sum z_je(c_j)}{u^2}-\frac{\sum z_j\widetilde{E}(c_j)}{u})\right),
\end{eqnarray*}
viewed as a morphism from
\[
\left((\rCn_*(A)\otimes\rCn_*(A^{\rm op}))\zz\uu\,\,,\,\, (b+\sum z_jb(c_j)+uB)\otimes1+1\otimes(b-\sum z_jb(c_j)+uB)\right)
\]
to $(\C\zz\uu\,\,,\,\,0)$, is 0-homotopic.
\end{lemma}
\noindent{\bf Proof.} In view of the previous lemma, the proof reduces to showing that
\[
\partial_u\cdot\STR\cdot\rCn_*(\rho)\cdot\ck-
\STR\cdot\rCn_*(\rho)\cdot\ck\cdot\left( {\conn}_{{\!\!{\partial_u}}}^{A}\otimes 1+ 1\otimes{\conn}_{{\!\!{\partial_u}}}^{A^{\rm op}}\right)
\]
is 0-homotopic. As we explained earlier (see Proposition \ref{kun} and the proof of Proposition \ref{thm1}), there exist odd operators
\[
H_1: \rCn_*(\End\,A)\uu\to\C\uu,\quad H_2: (\rCn_*(A)\otimes\rCn_*(A^{\rm op}))\uu\to\rCn_*(A\otimes  A^{\rm op})\uu
\]
such that $\partial_u\cdot\STR-\STR\cdot {\conn}_{{\!\!{\partial_u}}}^{\End\,A}=H_1\cdot (b+uB)$ and
\begin{eqnarray*}
&&{\conn}_{{\!\!{\partial_u}}}^{A\otimes  A^{\rm op}}\cdot(\ck+u\cK)-(\ck+u\cK)\cdot\left( {\conn}_{{\!\!{\partial_u}}}^{A}\otimes 1+ 1\otimes{\conn}_{{\!\!{\partial_u}}}^{A^{\rm op}}\right)=\\&&=(b+uB)\cdot H_2+H_2\cdot ((b+uB)\otimes1+1\otimes(b+uB)).
\end{eqnarray*}
Then
\begin{eqnarray*}
\partial_u\cdot\STR\cdot\rCn_*(\rho)\cdot\ck-
\STR\cdot\rCn_*(\rho)\cdot\ck\cdot\left( {\conn}_{{\!\!{\partial_u}}}^{A}\otimes 1+ 1\otimes{\conn}_{{\!\!{\partial_u}}}^{A^{\rm op}}\right)=\\
=\partial_u\cdot\STR\cdot\rCn_*(\rho)\cdot(\ck+u\cK)-
\STR\cdot\rCn_*(\rho)\cdot(\ck+u\cK)\cdot\left( {\conn}_{{\!\!{\partial_u}}}^{A}\otimes 1+ 1\otimes{\conn}_{{\!\!{\partial_u}}}^{A^{\rm op}}\right)=\\
=\partial_u\cdot\STR\cdot\rCn_*(\rho)\cdot(\ck+u\cK)-\STR\cdot {\conn}_{{\!\!{\partial_u}}}^{\End\,A}\cdot\rCn_*(\rho)\cdot(\ck+u\cK)+\\
+\STR\cdot\rCn_*(\rho)\cdot{\conn}_{{\!\!{\partial_u}}}^{A\otimes  A^{\rm op}}\cdot(\ck+u\cK)-\STR\cdot\rCn_*(\rho)\cdot(\ck+u\cK)\cdot\left( {\conn}_{{\!\!{\partial_u}}}^{A}\otimes 1+ 1\otimes{\conn}_{{\!\!{\partial_u}}}^{A^{\rm op}}\right)=\\
=H_1\cdot (b+uB)\cdot\rCn_*(\rho)\cdot(\ck+u\cK)+\STR\cdot\rCn_*(\rho)\cdot H_2\cdot ((b+uB)\otimes1+1\otimes(b+uB))=\\
=(H_1\cdot\rCn_*(\rho)\cdot(\ck+u\cK)+\STR\cdot\rCn_*(\rho)\cdot H_2)((b+uB)\otimes1+1\otimes(b+uB)).
\end{eqnarray*}
Therefore, it would suffice to show that
\begin{eqnarray*}
(H_1\cdot\rCn_*(\rho)\cdot(\ck+u\cK)+\STR\cdot\rCn_*(\rho)\cdot H_2)(\sum z_jb(c_j)\otimes1-1\otimes \sum z_jb(c_j))=0.
\end{eqnarray*}
We will prove that for any central closed even $c$
\begin{eqnarray}\label{h1}
\rCn_*(\rho)\cdot(\ck+u\cK)\cdot (b(c)\otimes1-1\otimes b(c))=0
\end{eqnarray}
(thus, an explicit form of $H_1$ is not needed) and 
\begin{eqnarray}\label{h2}
\STR\cdot\rCn_*(\rho)\cdot H_2\cdot(b(c)\otimes1-1\otimes b(c))=0.
\end{eqnarray}

To prove (\ref{h1}), we observe that
\begin{eqnarray}\label{shb}
\ck\cdot (b(c)\otimes1)=b(c\otimes1)\cdot\ck,\quad \ck\cdot (1\otimes b(c))=b(1\otimes c)\cdot\ck,
\end{eqnarray}
\begin{eqnarray}\label{Shb}
\cK\cdot (b(c)\otimes1)=b(c\otimes1)\cdot\cK,\quad \cK\cdot (1\otimes b(c))=b(1\otimes c)\cdot\cK.
\end{eqnarray}
We omit here the (quite elementary) proof of this claim; roughly, the main idea is that for any two sequences $x_1,\ldots, x_l$ and $y_1,\ldots, y_m$ there are obvious one-to-one correspondences between the following three sets: the set of sequences obtained from the shuffles of $x_1|\ldots| x_l|y_1|\ldots |y_m$ by inserting $z$ at a random position, the set of sequences obtained by shuffling all the sequences $x_1|\ldots|x_{i-1}|z|x_i|\ldots| x_l|y_1|\ldots |y_m$ ($i=0,\ldots,l+1$), and the set of sequences obtained by shuffling all the sequences $x_1|\ldots| x_l|y_1|\ldots|y_{i-1}|z|y_i|\ldots |y_m$ ($i=0,\ldots,m+1$).

Equalities (\ref{shb}) and (\ref{Shb}) yield
\begin{eqnarray*}
\ck\cdot (b(c)\otimes1-1\otimes b(c))=b(c\otimes1-1\otimes c)\cdot\ck,\\ \cK\cdot (b(c)\otimes1-1\otimes b(c))=b(c\otimes1-1\otimes c)\cdot\cK.
\end{eqnarray*}
Then
\begin{eqnarray*}
\rCn_*(\rho)\cdot(\ck+u\cK)\cdot (b(c)\otimes1-1\otimes b(c))=\rCn_*(\rho)\cdot b(c\otimes1-1\otimes c)\cdot(\ck+u\cK)\\=b(L_c-R_c)\cdot\rCn_*(\rho)\cdot (\ck+u\cK)=0.
\end{eqnarray*}

To prove (\ref{h2}), we need an explicit formula for $H_2$ obtained in \cite{Shk3}:
\[
H_2=\frac1{2u} H-\frac1{2u^2}(\ck+u\cK)\cdot((\U+u\V)\otimes1)\cdot(1\otimes B)
\]
where $H$ is an operator that produces tensors of the form $1\otimes1[\ldots]\in\rCn_*(A\otimes  A^{\rm op})$ and, therefore, does not contribute to $\STR\cdot\rCn_*(\rho)\cdot H_2\cdot(b(c)\otimes1-1\otimes b(c))$ (see the beginning of the proof of Proposition \ref{thm1}). Finally,
\[
\STR\cdot\rCn_*(\rho)\cdot (\ck+u\cK)\cdot((\U+u\V)\otimes1)\cdot(1\otimes B)\cdot(b(c)\otimes1-1\otimes b(c))=0
\]
by Proposition \ref{le1}, (\ref{ii}), and (\ref{h1}). 
\hfill $\blacksquare$

\end{document}